\documentclass[10pt,a4paper]{smfart}
\usepackage[latin1]{inputenc}
\usepackage{amsmath}
\usepackage{pstricks}
\usepackage{amsfonts}
\usepackage{amssymb}
\usepackage{amsthm}
\usepackage{lmodern}
\usepackage[all]{xy}
\usepackage{mathrsfs}
\usepackage{tikz}
\usepackage[english,francais]{babel}
\newtheorem{defi}{Definition}[section]
\newtheorem{thm}[defi]{Theorem}

\newtheorem{coro}[defi]{Corollary}

\newtheorem{rem}[defi]{Remark}

\author{Jean-Christophe SAN SATURNINO}
\title[Ramification theory and key polynomials]{Ramification theory and key polynomials}

\address{Universit\'e Toulouse III Paul Sabatier\\
Institut de Math\'ematiques de Toulouse\\ 
118 route de Narbonne\\
31062 Toulouse cedex 9 (France)}
\email{san@math.ups-tlse.fr}
\urladdr{http://www.math.univ-toulouse.fr/$\sim$san/}

\begin{document}
\selectlanguage{english}
\maketitle
\begin{abstract}
For a simple, normal and finite extension of a valued field, we prove that we can related the order of the ramification group of the field extension and the set of key polynomials associated to the extension of the valuation. More precisely, the order of this group can be expressed in terms of a product of a power of the characteristic of the residue field of the valuation and the effective degrees of the key polynomials. We also give a condition on the order of the ramification group so that there is no limit key polynomials for a valuation of rank one. This condition also allow us to have a monomialization theorem.
\end{abstract}
\setcounter{section}{-1}
\section{Introduction}
\indent Since the work of S. Abhyankar, F.-V. Kuhlmann, V. Cossart and O. Piltant, we know the importance to studying the decomposition, ramification and inertia groups of a valuation in the problem of local uniformization. In \cite{jccar0}, we saw the importance of the key polynomials to obtain a local uniformization for a valuation of rank one. If the characteristic of the residual field is $0$ then we obtain the result, otherwise it is sufisent to monomialize the first limit key polynomial wich is of the form:
$$Q_\omega=X^{p^e}+\sum\limits_{i=0}^{e-1}c_{p^i}X^{p^i}+c_0.$$
\\ \\\indent In \cite{jcdefaut}, we related the defect of a simple and finite extension of valued fields with the effective degree of the key polynomials. For example, if the valuation is of rank $1$ and there is only one limit key polynomial with the others of degree$1$, the defect is exactly the degree of the limit key polynomial.
\\\indent This paper complete \cite{jcdefaut} in our program to connect the the ramification theory, the defect of an finite extension and the works of F.-V. Kuhlmann in local uniformization with the key polynomials. For a simple, normal and finite extension of a valued field, we prove that we can related the order of the ramification group of the field extension and the set of key polynomials associated to the extension of the valuation.
\\ \\\indent In the first section we recall the definitions of the  decomposition, ramification and inertia groups. We also give some relations between the dimension of the different fields and the order of the groups of valuation.
\\\indent In the second chapter we recall the results of \cite{jcdefaut} and interpret them in terms of the order of ramification group. In particular we give a condition on the order of the ramification group to have no defect or no limit key polynomials.
\\\indent In the last section we give some conditions to obtain a local uniformization theorem with asumptions on the group of ramification. 
\\ \\\noindent\textbf{Notation.} Let $\nu$ be a valuation of a field $K$. We write $R_\nu=\lbrace f\in K\:\vert\:\nu(f)\geqslant 0\rbrace$, this is a local ring whose maximal ideal is $m_\nu=\lbrace f\in K\:\vert\:\nu(f)> 0\rbrace$. We then denote by $k_\nu=R_\nu/m_\nu$ the residue field of $R_\nu$ and $\Gamma_\nu=\nu(K^*)$.
\\For a field $K$, we will denote by $\overline K$ an algebraic closure of $K$. For an algebraic extension $L\vert K$ we will denote by $(L\vert K)^{sep}$ a separable closure of $K$ in $L$ (or more simply by $K^{sep}$ if no confusion is possible) and by $Aut(L\vert K)$ 
 the group of automorphisms of $L\vert K$ (or $Gal(L\vert K)$ if $L\vert K$ is a Galois extension). \\For a finite extension $L\vert K$, we denote by $[L:K]_{sep}$ the \textit{separable degree} of $L\vert K$ and by $[L:K]_{ins}=[L:K]/[L:K]_{sep}$ the \textit{inseparable degree} of $L\vert K$; this a power of $car(K)$.
\\If $R$ is a ring and $I$ an ideal of $R$, we will denote by $\widehat{R}^I$ the \textit{ $I$-adic completion} of $R$. When $(R,\mathfrak m)$ is a local ring, we will say the \textit{completion} of $R$ instead of the $\mathfrak m$-adic completion  of $R$ and we will denote it by $\widehat{R}$.
\\For all $P\in Spec(R)$, we note $\kappa(P)=R_P/PR_P$ the residue field of $R_P$.
\\For $\alpha\in\mathbb Z^n$ and $u=(u_1,...,u_n)$ a $n$-uplet of elements of $R$, we write:
\[u^\alpha=u_1^{\alpha_1}...u_n^{\alpha_n}.\]
For $P,Q\in R\left[ X\right] $ with $P=\sum\limits_{i=0}^{n}a_{i}Q^{i}$ and $a_{i}\in R[X]$ such that the degree of $a_{i}$ is strictly less than $Q$, we write:
\[d_{Q}^{\:\circ}(P)=n.\]
If $Q=X$, we will simply write $d^{\:\circ}(P)$ instead of $d_{X}^{\:\circ}(P)$.
\\Finally, if $R$ is a domain, we denote by $Frac(R)$ its quotient field.
\section{Decomposition, inertia and ramification groups}
\indent We follow the definitions of \cite{kulhlivre}, Chapter 7.
\begin{defi}
Let $(K,\nu)\hookrightarrow (L,\mu)$ be a normal algebraic extension of valued fields. 
\begin{enumerate}
\item The \textbf{decomposition group of} $\boldsymbol{L\vert K}$ is:
\[G^d(L\vert K,\mu)=\{\sigma\in Aut(L\vert K)\:\vert\:\forall\alpha\in L,\:\mu(\sigma (\alpha))=\mu(\alpha)\}.\]
\item The \textbf{inertia group of} $\boldsymbol{L\vert K}$ is:
\[G^i(L\vert K,\mu)=\{\sigma\in Aut(L\vert K)\:\vert\:\forall\alpha\in R_\mu,\:\mu(\sigma (\alpha)-\alpha)>0\}.\]
\item The \textbf{ramification group of} $\boldsymbol{L\vert K}$ is:
\[G^r(L\vert K,\mu)=\{\sigma\in Aut(L\vert K)\:\vert\:\forall\alpha\in R_\mu,\:\mu(\sigma (\alpha)-\alpha)>\mu(\alpha)\}.\]
\end{enumerate}
When no confusion is possible, we will denote respctiveley by $G^d$, $G^i$ and $G^r$ the groups $G^d(L\vert K,\mu)$, $G^i(L\vert K,\mu)$ and $G^r(L\vert K,\mu)$.
\end{defi}
\begin{rem}\label{normgroup}
\textup{$G^r(L\vert K,\mu)\lhd G^i(L\vert K,\mu)\lhd G^d(L\vert K,\mu)\lhd Aut(L\vert K)$.}
\end{rem}
\begin{defi}
Let $(K,\nu)\hookrightarrow (L,\mu)$ be a normal algebraic extension of valued fields. 
\begin{enumerate}
\item The fixed field of $G^d(L\vert K,\mu)$ in the separable closure of $K$ in $L$ is called the \textbf{decomposition field of} $\boldsymbol{K}$ \textbf{in} $\boldsymbol L$ and will be denote by $K^d$. We also write $\mu^d=\mu_{\vert K^d}$.
\item The fixed field of $G^i(L\vert K,\mu)$ in the separable closure of $K$ in $L$ is called the \textbf{inertia field of} $\boldsymbol{K}$ \textbf{in} $\boldsymbol L$ and will be denote by $K^i$. We also write $\mu^i=\mu_{\vert K^i}$.
\item The fixed field of $G^r(L\vert K,\mu)$ in the separable closure of $K$ in $L$ is called the \textbf{ramification field of} $\boldsymbol{K}$ \textbf{in} $\boldsymbol L$ and will be denote by $K^r$. We also write $\mu^r=\mu_{\vert K^r}$.
\end{enumerate}
\end{defi}
\begin{rem}
\textup{The ramification, inertia and decomposition fields are separable over $K$. By the remark \ref{normgroup}, the extensions $K^i\vert K^d$ and $K^r\vert K^d$ are Galois extensions.
\\If $L=\overline K$ then $K^d$ is an henselization of $K$ denote by $K^h$.}
\end{rem}
\begin{thm}\label{thmgal}
Let $(K,\nu)\hookrightarrow (L,\mu)$ be a normal and finite extension of valued fields. Write $e=[\Gamma_\mu:\Gamma_\nu]$, $f=[k_\mu:k_\nu]$, $p=char(k_\nu)$ and $G=Aut(L\vert K)$. Let $g$ be the number of distinct extensions of $\nu$ from $K$ to $L$ and $d=d_{L\vert K}(\mu,\nu)$ be the defect of the extension $L\vert K$ in $\mu$ (see Definition 2.3 of \cite{jcdefaut}).
\begin{enumerate}
\item $G=Gal(K^{sep}\vert K)$, $G^d=Aut(L\vert K^d)=Gal(K^{sep}\vert K^d)$, $G^i=Aut(L\vert K^i)=Gal(K^{sep}\vert K^i)$, $G^r=Aut(L\vert K^r)=Gal(K^{sep}\vert K^r)$.
\item $[L:K]=g\times d\times e\times f$.
\item $K^d\vert K$ is an immediate extension and $g=[G:G^d]=[K^d:K]$.
\item $K^i\vert K^d$ is a Galois extension, $[\Gamma_{\mu^i}:\Gamma_{\mu^d}]=1$, 
$[k_{\mu^i}:k_{\mu^d}]=[K^i:K^d]=f_0$.
\item $K^r\vert K^i$ is an abelian extension, $[\Gamma_{\mu^r}:\Gamma_{\mu^i}]=[K^r:K^i]=e_0$, $[k_{\mu^r}:k_{\mu^i}]=1$.
\item $K^{sep}\vert K^r$ is a $p$-extension, $[\Gamma_{\mu}:\Gamma_{\mu^r}]=p^t$, $[k_{\mu}:k_{\mu^r}]=p^s$ with $t,s\geqslant 0$ and $[K^{sep}:K^r]=\vert G^r\vert=p^u$.
\end{enumerate}
\end{thm}
\noindent\textit{Proof}: For a detailed proof, one can consult \cite{kulhlivre}, Chapter 7. We will give some ideas of proof.
\begin{enumerate}
\item It is clear because $L\vert K^{sep}$ is a purely inseparable extension and by using Galois theory.
\item By Lemma 7.46 of \cite{kulhlivre}, we know that:
\[\left[L:K\right]=\sum\limits_{i=1}^gd_ie_if_i,\]
where $d_i=d_{L\vert K}(\mu_i,\nu)$, $e_i=\left[\Gamma_{\mu_i}:\Gamma_\nu\right]$, $f_i=\left[k_{\mu_i}:k_\nu\right]$ and $\mu_i$ all the extensions of $\nu$ to $L$, $i\in\{1,...,g\}$. The extension $L\vert K$ is a normal extension so, for $i\neq j$, $e_i=e_j$, $f_i=f_j$ and $d_i=d_j$. Then we deduce the equality.
\item To prove that $g=[G:G^d]$, it is suffice to observe that, for $\sigma,\sigma'\in G$,  $\mu(\sigma(\alpha))=\mu(\sigma'(\alpha))$ if and only if $\sigma\sigma'^{-1}\in G^d$. The second equality comme from 1. and Galois theory.
\\To show that $K^d\vert K$ is an immediate extension, we need to prove that: $$[k_{\mu^d}:k_{\nu}]=1=[\Gamma_{\mu^d}:\Gamma_{\nu}].$$ The first equality come from the fact that $\overline\alpha=\overline{Tr_{K^d\vert K}(\alpha)}$ for $\overline\alpha\in k_{\mu^d}\setminus\{0\}$. To prove the second equality, take $\gamma\in\Gamma_{\mu^d}$ and construct an element of $K^d$ such that is minimal polynomial over $K$ have only one root of value $\gamma$: this element.
\item $K^i\vert K^d$ is a Galois extension by Galois theory because  $G^i(L\vert K,\mu)\lhd G^d(L\vert K,\mu)$ and we have $Gal(K^i\vert K^d)\simeq G^d/G^i$. Using the separable closure of $k_{\mu^d}$ in $k_\mu$ (which is $k_{\mu^i}$), we can prove that $G^d\rightarrow Aut(k_\mu\vert k_\nu)$ is surjective. So, by definition of $G^i$, $Gal(K^i\vert K^d)\simeq G^d/G^i\simeq Gal(k_\mu\vert k_{\mu^d})=Aut(k_\mu\vert k_\nu)$. In the same way, we can show that $Gal(k_\mu\vert k_{\mu^d})\simeq Gal(k_{\mu^i}\vert k_{\mu^d})$. Finally from the fundamental inequality, we obtain:
$$[K^i:K^d]=\left\vert Gal(K^i\vert K^d)\right\vert=\left\vert Gal(k_{\mu^i}\vert k_{\mu^d})\right\vert\leqslant [k_{\mu^i}: k_{\mu^d}]\leqslant [K^i:K^d].$$
\item$K^r\vert K^i$ is a Galois extension by Galois theory because  $G^r(L\vert K,\mu)\lhd G^i(L\vert K,\mu)$ and we have $Gal(K^r\vert K^i)\simeq G^i/G^r$. Assuming that 6. is true, we can show that there exists an embedding between $Gal(K^r\vert K^i)$ and $Hom((\Gamma_{\mu^r}/\Gamma_{\mu^i})_{p'},k_{\mu^r}^\times)$ where $\Gamma_{\mu^r}/\Gamma_{\mu^i}=(\Gamma_{\mu^r}/\Gamma_{\mu^i})_{p}\oplus (\Gamma_{\mu^r}/\Gamma_{\mu^i})_{p'}$ with $(\Gamma_{\mu^r}/\Gamma_{\mu^i})_{p}$ a $p$-group and $(\Gamma_{\mu^r}/\Gamma_{\mu^i})_{p'}$ a torsion group wich the orders of all elements are prime to $p$. We deduce from this that $K^r\vert K^i$ is an abelian extension. Using the fundamental inequality, we obtain:
\begin{align*}
[K^r:K^i]=\vert Gal(K^r\vert K^i)\vert&\leqslant\vert Hom((\Gamma_{\mu^r}/\Gamma_{\mu^i})_{p'},k_{\mu^r}^\times)\vert\\&\leqslant \left\vert Hom\left((\Gamma_{\mu^r}/\Gamma_{\mu^i})_{p'},\overline{k_{\mu^r}}^\times\right)\right\vert=\vert (\Gamma_{\mu^r}/\Gamma_{\mu^i})_{p'}\vert\\&\leqslant\vert\Gamma_{\mu^r}/\Gamma_{\mu^i}\vert\\&\leqslant [K^r:K^i].
\end{align*}
We have proved that $[\Gamma_{\mu^r}:\Gamma_{\mu^i}]=[K^r:K^i]$ and using the fundamental inequality, we obtain that $[k_{\mu^r}:k_{\mu^i}]=1$.
\item By Lemma 7.15 of \cite{kulhlivre}, it is sufficient to show that $G^r$ is a $p$-group, and we will have $[\Gamma_{\mu}:\Gamma_{\mu^r}]=p^t$. Since $k_{\mu^r}=k_{\mu^i}=k_{\nu}^{sep}$, then $[k_{\mu}:k_{\mu^r}]=[k_{\mu}:k_{\nu}]_{ins}$ is a power of $p$. To show that $G^r$ is a $p$-group, proceed by contradiction. Take an element of order a prime $q\neq p$. Consider the fixed field of this element, then $K^r$ is a cyclic extension of this field, of degree $q$. Let $x$ a primitive element, we can suppose that this trace is $0$ because $q\neq 0$. On the other hand, the sum of the residue of all the conjugates of $x$ divide by $x$ is $q$, so it is a contradiction with we consider his trace.
\end{enumerate}\qed
\\ \\As in \cite{kulhlivre}, Chapter 7, we can summarize the Theorem \ref{thmgal} in this table with $n=[L:K]$:
\newline
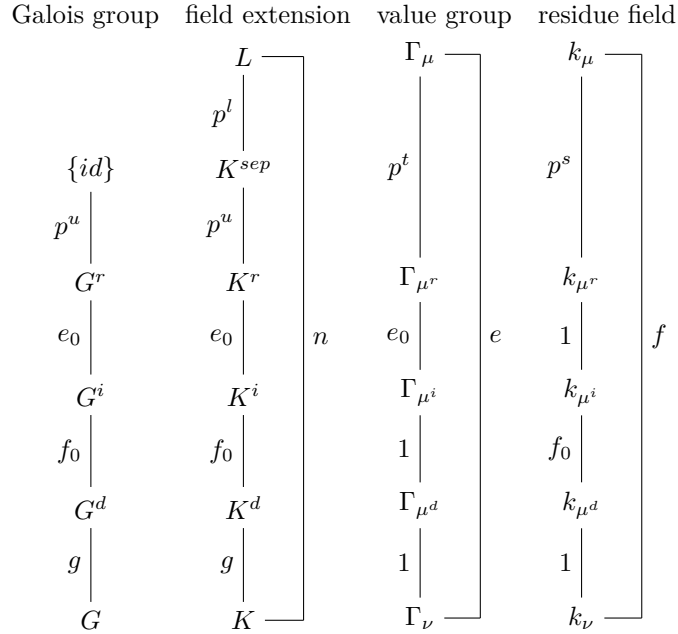
\begin{center}
\begin{tabular}{cccc}

Galois group&field extension& value group & residue field\\
\begin{tikzpicture}[node distance = 1.5cm, auto]
      \node (G) {$G$};
      \node (Gd) [above of=G] {$G^d$};
      \node (Gi) [above of=Gd] {$G^i$};
      \node (Gr) [above of=Gi] {$G^r$};
      \node (id) [above of=Gr] {$\{id\}$};
      \draw[-] (G) to node {$g$} (Gd);
      \draw[-] (Gd) to node {$f_0$} (Gi);
      \draw[-] (Gi) to node {$e_0$} (Gr);
      \draw[-] (Gr) to node {$p^u$} (id);
    
      \end{tikzpicture}
      &
      \begin{tikzpicture}[node distance = 1.5cm, auto]
      \node (K) {$K$};
      \node (Kd) [,above of=K] {$K^d$};
      \node (Ki) [above of=Kd] {$K^i$};
      \node (Kr) [above of=Ki] {$K^r$};
      \node (Ksep) [above of=Kr] {$K^{sep}$};
      \node (L) [above of=Ksep] {$L$};
      \draw[-] (K) to node {$g$} (Kd);
      \draw[-] (Kd) to node {$f_0$} (Ki);
      \draw[-] (Ki) to node {$e_0$} (Kr);
      \draw[-] (Kr) to node {$p^u$} (Ksep);
      \draw[-] (Ksep) to node {$p^l$} (L);
      \draw (K) -- (0.8,0); 
       \draw[-] (0.8,0)  to node[swap] {$n$} (0.8,7.5); 
      \draw (0.8,7.5) -- (L); 
      \end{tikzpicture}
      &\begin{tikzpicture}[node distance = 1.5cm, auto]
      \node (g) {$\Gamma_{\nu}$};
      \node (gd) [above of=g] {$\Gamma_{\mu^d}$};
      \node (gi) [above of=gd] {$\Gamma_{\mu^i}$};
      \node (gr) [above of=gi] {$\Gamma_{\mu^r}$};
      \node (gsep) [above of=gr] {};
      \node (gmu) [above of=gsep] {$\Gamma_{\mu}$};
      \draw[-] (g) to node {$1$} (gd);
      \draw[-] (gd) to node {$1$} (gi);
      \draw[-] (gi) to node {$e_0$} (gr);
       \draw[-] (gr) to node {$p^t$} (gmu);
        \draw (g) -- (0.8,0); 
       \draw[-] (0.8,0)  to node[swap] {$e$} (0.8,7.5); 
      \draw (0.8,7.5) -- (gmu); 
      \end{tikzpicture}
      &\begin{tikzpicture}[node distance = 1.5cm, auto]
      \node (g) {$k_{\nu}$};
      \node (gd) [above of=g] {$k_{\mu^d}$};
      \node (gi) [above of=gd] {$k_{\mu^i}$};
      \node (gr) [above of=gi] {$k_{\mu^r}$};
      \node (gsep) [above of=gr] {};
      \node (gmu) [above of=gsep] {$k_{\mu}$};
      \draw[-] (g) to node {$1$} (gd);
      \draw[-] (gd) to node {$f_0$} (gi);
      \draw[-] (gi) to node {$1$} (gr);
       \draw[-] (gr) to node {$p^s$} (gmu);
      \draw (g) -- (0.8,0); 
       \draw[-] (0.8,0)  to node[swap] {$f$} (0.8,7.5); 
      \draw (0.8,7.5) -- (gmu); 
      \end{tikzpicture}

\end{tabular}
\end{center}
~
\begin{coro}\label{corogaldef}
Under the same assumptions and notations of the Theorem \ref{thmgal}, we have:
\begin{enumerate}
\item $e=e_0\times p^t$ with $p\not\vert\; e_0$ and $f=f_0\times p^s$.
\item $d=p^{u+l-s-t}$ where $p^l=[L:K]_{ins}$.
\end{enumerate}
\end{coro}
\begin{rem}
\textup{If $L\vert K$ is a finite Galois extension, then $l=0$ an the Corollary \ref{corogaldef} is the same as the Corollary of Theorem 25, Ch. VI of \cite{zarsam2}. }
\end{rem}
\noindent\textit{Proof}: $e=e_0\times p^t$ and $f=f_0\times p^s$ come from the precedent table and $p\not\vert\; e_0$ because $G^r$ is the unique $p$-Sylow of $G^i$ and $\vert G^i\vert=e_0\times p^u$. Finally by 1. of Theorem \ref{thmgal} and the precedent table, since:
\[g\times d\times e\times f=[L:K]=p^{l+u}\times e_0\times f_0\times g,\]
then: 
$$d=p^{u+l-s-t}.$$
\qed

\section{Key polynomials and ramification group}
\indent For a basic definition of key polynomials see Definition 3.1 of \cite{jcdefaut}, for more details see \cite{spivamahboub}. Here, we suppose known the theory of key polynomials, we only recall the link with the defect as in \cite{jcdefaut}.
\\\indent Let $\lbrace Q_{l}\rbrace_{l\in\Lambda}$ be a complete set of key polynomials, for $l\in\Lambda$ having a predecessor, write:
$\alpha_{l}=d_{Q_{l-1}}^{\:\circ}(Q_{l}).$
If $l=\omega n$ is a limit ordinal, $n\in\mathbb N^*$, denote by $\alpha_{l}=d_{Q_{l_0}}^{\:\circ}(Q_{l})$ where $l_0=\min\lbrace m\geqslant 1\:\vert\:\alpha_{\omega(n-1)+m}=1\rbrace$.
\begin{defi}
Let $K$ be a field and $\mu$ a valuation of $K[x]$. For $h\in K[x]$, consider its $i$-standard expansion $h=\sum\limits_{j=0}^{s_{i}}c_{j,i}Q_{i}^j$. We call by the $\boldsymbol{i}$\textbf{-th effective degree of h} the natural number:
\[\delta_i(h)=\max \lbrace j\in\lbrace 0,...,s_{i}\rbrace\:\vert\:j\beta_i+\mu\left(c_{j,i}\right)=\mu_i(h)\rbrace,\]
where:
\[\mu_{i}(h)=\min_{0\leqslant j\leqslant s_{i}}\lbrace j\mu(Q_{i})+\mu(c_{j,i})\rbrace.\]
By convention, $\delta_i(0)=-\infty$.
\end{defi}
\begin{rem}
\textup{Remind that, by Proposition 5.2 of \cite{spivamahboub}, for $l\in\Lambda$ an ordinal number, the sequence $(\delta_{l+i}(h))_{i\in\mathbb N^*}$ decreases. Thus there exists $i_0\in\mathbb N^*$ such that $\delta_{l+i_0}(h)=\delta_{l+i_0+i}(h)$, for all $i\geqslant 1$ and we denote this common value by $\delta_{l+\omega}(h)$. From here until the end, we write $\delta_{l+\omega}=\delta_{l+\omega}(Q_{l+\omega})$.}
\end{rem}
\begin{thm}(\cite{jcdefaut}, Corollary 4.5)\label{defautdelta}
Let $(K,\nu)$ be a valued field and $L$ be a finite and simple extension of $K$. Write $\mu^{(1)},...,\mu^{(g)}$ the different extensions of $\nu$ on $L$, its corresponds to a (pseudo-)valuation of $K\left[ x\right]$ denoted by the same way. Consider $\lbrace Q_l^{(i)} \rbrace_{l\in\Lambda^{(i)}}$ the set of key polynomials associated to $\mu^{(i)}$ and $n_0^{(i)}\in\mathbb N^*$ the smallest possible such that $\Lambda^{(i)}\leqslant \omega n_0^{(i)}$, $1\leqslant i\leqslant g$. Then:
\[d_{L\vert K}(\mu^{(i)},\nu)=\prod\limits_{j=1}^{n_0^{(i)}}d_{\omega j}^{(i)}.\]
We deduced that:
\[ \left[ L:K\right]=\sum\limits_{i=1}^ge_if_id_{\omega}^{(i)}d_{\omega 2}^{(i)}...d_{\omega n_0^{(i)}}^{(i)},\]
where $e_i=\left[\Gamma_{\mu^{(i)}}:\Gamma_\nu\right]$, $f_i=\left[k_{\mu^{(i)}}:k_\nu\right]$, $d_{\omega j}^{(i)}=\delta_{\omega j}^{(i)}$ for $j<n_0^{(i)}$ and:
\[d_{\omega n_0^{(i)}}^{(i)}=\left \{ \begin{array}{ccl}  \delta_{\omega n_0^{(i)}}^{(i)} & \textup{si} & \Lambda=\omega n_0^{(i)}\textup{ et } \sharp\lbrace m\geqslant 1\:\vert\:\alpha_{\omega(n_0^{(i)}-1)+m}^{(i)}=1\rbrace=+\infty \\  1  & \textup{si} & \Lambda<\omega n_0^{(i)} \textup{ ou } \Lambda=\omega n_0^{(i)}\textup{ et } \sharp\lbrace m\geqslant 1\:\vert\:\alpha_{\omega(n_0^{(i)}-1)+m}^{(i)}=1\rbrace<+\infty \end{array} \right.\] 
\end{thm}
\begin{rem}
\textup{Denote by $p=char(k_\nu)$. By a suggestion of G. Leloup, since $d_{\omega j}^{(i)}\geqslant p$ for $j\in\{1,...n_0^{(i)}-1\}$ and $d_{\omega n_0^{(i)}}^{(i)}\geqslant 1$, we have:
\[d_{L\vert K}(\mu^{(i)},\nu)\geqslant p^{n_0^{(i)}-1}.\]
We deduce that:
\[n_0^{(i)}\leqslant\log_p\left(d_{L\vert K}(\mu^{(i)},\nu)\right)+1.\]
This result brings more precision than the inequality given in \cite{spivaherrera}.}
\end{rem}
From now until the end of this section we will use the notations of the previous section.
\begin{thm}\label{thmprinc}
Consider the same assumptions as Theorem \ref{defautdelta} and assume more that $L\vert K$ is a normal extension. Then, for all $i\in\{1,...,g\}$, we have:
\[\vert G^r\vert = p^{s+t-l}\times\prod\limits_{j=1}^{n_0^{(i)}}d_{\omega j}^{(i)}.\]
\end{thm}
\noindent\textit{Proof}: By 2. of Corollary \ref{corogaldef}:
$$\vert G^r\vert=p^{s+t-l}\times d_{L\vert K}(\mu^{(i)},\nu).$$ 
To conclude it is sufficient to apply Theorem \ref{defautdelta}.\\\qed
\begin{coro}
With the same assumptions of Theorem \ref{thmprinc}, $d_{L\vert K}(\mu^{(i)},\nu)=1$, for all $i\in\{1...,g\}$, if and only if $\vert G^r\vert = p^{s+t-l}$.
\end{coro}
\indent In some cases, we can express the order of the ramification group only in terms of key polynomials, essentially with the degree of the key polynomials and the effective degree.
\begin{coro}
Consider the same assumptions as Theorem \ref{defautdelta} and assume more that $L\vert K$ is a Galois extension. Assume that $K^r=K^d$ or equivalently $G^r=G^i=G^d$. Then, for all $i\in\{1,...,g\}$, there exists an index $i_0\in\Lambda^{(i)}$ having a predecessor, such that:
\[\vert G^r\vert = d^{\:\circ}\left(Q_{\omega(n_0^{(i)}-1)+i_0}^{(i)}\right)\times d_{\omega n_0^{(i)}}^{(i)}.\]
\end{coro}
\noindent\textit{Proof}: The extension $L\vert K$ is Galois then $l=0$. The assumption $K^r=K^d$ is equivalent to $e_0=f_0=1$. Thus $e=p^t$, $f=p^s$ and, by Theorem \ref{thmprinc}: $$\vert G^r\vert=\left(e\times f\times \prod\limits_{j=1}^{n_0^{(i)}-1}d_{\omega j}^{(i)}\right)d_{\omega n_0^{(i)}}^{(i)}.$$ But in the proof of Corollary 4.3 of \cite{jcdefaut}, we have seen that, as a consequence of Proposition 2.9 of \cite{vaquie3}, there exists an index $i_0\in\Lambda^{(i)}$ having a predecessor, such that:
$$d^{\:\circ}\left(Q_{\omega(n_0^{(i)}-1)+i_0}^{(i)}\right)=e\times f\times\prod\limits_{j=1}^{n_0^{(i)}-1}d_{\omega j}^{(i)}.$$\qed
\begin{thm}\label{pasdepolycle}
Consider the same assumptions as Theorem \ref{thmprinc} and assume more that $rk(\nu)=1$. If $\vert G^r\vert = p^{s+t-l}$ then the set of key polynomials $\lbrace Q_l^{(i)} \rbrace_{l\in\Lambda^{(i)}}$ associated to $\mu^{(i)}$ have no limit key polynomials, ie: $\Lambda^{(i)}\subseteq\mathbb N^*$.
\end{thm}
\noindent\textit{Proof}: If $\vert G^r\vert = p^{s+t-l}$ then, by definition, $n_0^{(i)}=1$ and $d_\omega^{(i)}=1$. We have three possibilities:
\begin{enumerate}
\item $\Lambda^{(i)}<\omega$, we have nothing to proove;
\item $\Lambda^{(i)}=\omega$ and $\sharp\lbrace m\geqslant 1\:\vert\:\alpha_{m}^{(i)}=1\rbrace<+\infty$, we conclude with the Proposition 3.19 of \cite{jcdefaut};
\item $\Lambda^{(i)}=\omega$ and $\sharp\lbrace m\geqslant 1\:\vert\:\alpha_{m}^{(i)}=1\rbrace=+\infty$, we conclude with the Proposition 3.20 and the Proposition 3.18 of \cite{jcdefaut} because $\delta_\omega^{(i)}=d_\omega^{(i)}=1$.
\end{enumerate}\qed
\begin{rem}
\textup{In the situation of the Theorem \ref{pasdepolycle}, the field is defectless and we can also apply directly the Proposition 5.1 of \cite{jcdefaut}.}
\end{rem}
\section{Ramification and local uniformization}
\indent Let $(R,\mathfrak{m},k)$ be a local complete regular equicharacteristic ring of dimension $n$ with $\mathfrak{m}=\left(u_1,...,u_{n}\right)$. Let $\nu$ be a valuation of $K=Frac(R)$, centered on $R$, of value group $\Gamma$ and $\Gamma_1$ the smallest non-zero isolated subgroup of $\Gamma$. Write: 
\[H=\lbrace f\in R\:\vert\: \nu(f)\notin\Gamma_{1}\rbrace.\]
$H$ is a prime ideal of $R$ (see the proof of Theorem 6.2 of \cite{jcdefaut}). Moreover suppose that:
\[n=e(R,\nu)=emb.dim\left(R/H\right),\]
 that is to say:
 \[H\subset\mathfrak{m}^2.\]
Write $r=r(R,u,\nu)=\dim_\mathbb Q\left(\sum\limits_{i=1}^n\mathbb Q\nu(u_i)\right)$.
\\ The valuation $\nu$ is unique if $ht(H)=1$; it is the composition of the valuation $\mu:L^{*}\rightarrow\Gamma_{1}$ of rank $1$ centered on $R/H$, where $L=Frac(R/H)$, with the valuation $\theta :K^{*}\rightarrow \Gamma / \Gamma_{1}$, centered on $R_{H}$, such that $k_{\theta}\simeq \kappa(H)$.
\\By abuse of notation, for $f\in R$, we will denote by $\mu(f)$ instead of $\mu(f\mod H)$.
By the Cohen's theorem, we can suppose that $R$ is of the form:
\[R=k\left[ \left[ u_{1},...,u_{n}\right] \right].\]
For $j\in \lbrace r+1,...,n\rbrace$, write $\lbrace Q_{j,i}\rbrace_{i\in\Lambda_{j}}$ the set of key polynomials of the extension $k\left( \left( u_{1},...,u_{j-1}\right) \right)\hookrightarrow k\left( \left( u_{1},...,u_{j-1}\right) \right)(u_{j})$, $\textbf{Q}_{j,i}=\left\lbrace Q_{j,i'}\vert i'\in\Lambda_{j},i'<i\right\rbrace $, $\Gamma^{(j)}$ the value group of $\nu_{\vert k\left( \left( u_{1},...,u_{j}\right) \right)}$ and $\nu_{j,i}$ the $i$-troncation of $\nu$ for this extension.
\\ \\\indent For the definition of local framed sequences, one may consult D\'efinition 7.1 and the sections 4.1 and 4.2 of \cite{jccar0}.
\begin{thm}\label{thmeclatformcar0}
Suppose that, for $R_{n-1}=k\left[ \left[ u_{1},...,u_{n-1}\right] \right]$ we have:
\begin{enumerate}
\item 
\begin{enumerate}
\item Or $H\cap R_{n-1}\neq (0)$ and there exists a local framed sequence $(R_{n-1},u)\rightarrow (R',u')$ such that: $$e(R',\nu)<e(R_{n-1},\nu);$$
\item Or $H\cap R_{n-1}=(0)$ and for all $f\in R_{n-1}$, there exists a local framed sequence $(R_{n-1},u)\rightarrow (R',u')$ such that $f$ is a monomial in $u'$ times a unit of $R'$.
\end{enumerate}
\item The local framed sequence $(R_{n-1},u)\rightarrow (R',u')$ of (1) can be chosen defined over $T$.
\end{enumerate}
Moreover suppose that the ramification group of the extension $k\left( \left( u_{1},...,u_{n-1}\right) \right)\hookrightarrow k\left( \left( u_{1},...,u_{n-1}\right) \right)\left[u_n\right]/H$ have order $p^{s+t-l}$ with $p^s=[k_{\mu}:k_{\mu^r}]$, $p^t=[\Gamma_{\mu}:\Gamma_{\mu^r}]$ and $p^l$ the inseparable degree of the extension. Then the assumptions 1. and 2. are true with $R$ instead of $R_{n-1}$.
\end{thm}
\noindent\textit{Proof}: The proof is the same as the proofs on Theorem 5.1 and 7.2 of \cite{jccar0}. With the assumptions of the Theorem \ref{thmeclatformcar0}, we can use the Proposition 5.2 of \cite{jccar0}: $H$ is generated by a irreducible monic polynomial in $u_n$. Since the order of the ramification group of the extension $k\left( \left( u_{1},...,u_{n-1}\right) \right)\hookrightarrow k\left( \left( u_{1},...,u_{n-1}\right) \right)\left[u_n\right]/H$ is $p^{s+t-l}$, by Theorem \ref{pasdepolycle}, the set of key polynomials $\lbrace Q_{j,i}\rbrace_{i\in\Lambda_{j}}$ has not limit key polynomial. To conclude it is sufficient to apply Theorem 7.2 of \cite{jccar0}.\\\qed
\begin{rem}
\textup{In \cite{jcthese} and \cite{jccar0} we saw that the problem of local uniformization is reduced to monomialize the first limit key polynomial $Q_\omega\in K[X]$. where $(K,\nu)$ is valuated field with $rk(\nu)=1$ and we have a local uniformization property on $K$. We know that we can suppose that:
$$Q_\omega=X^{p^e}+\sum\limits_{i=0}^{e-1}c_{p^i}X^{p^i}+c_0,$$
where $p^e$ is equal to the defect. In this situation, if we write $L=K[X]/(Q_\omega)$, we have for the extension $L\vert K$, $e=f=1$. This extension have defect and if we take a Galois closure of $L$, we obtain that the order of $G^r$ is exactely the defect. So, if we want to investigate a way to obtain a local uniformization theorem studying the ramification group of an extension, this the only situation which is problematic.}
\end{rem}

\bibliographystyle{plain}
\bibliography{biblio2}

\begin{thebibliography}{1}

\bibitem{spivamahboub}
F.~J. Herrera~Govantes, W.~Mahboub, M.~A. Olalla~Acosta, and M.~Spivakovsky.
\newblock \textit{Key polynomials for simple extensions of valued fields},
  2014.
\newblock preprint math.AG/arXiv:1406.0657.

\bibitem{spivaherrera}
F.~J. Herrera~Govantes, M.~A. Olalla~Acosta, and M.~Spivakovsky.
\newblock \textit{Valuations in algebraic field extensions}.
\newblock {\em \textup{J. Algebra}}, 312(2):1033--1074, 2007.

\bibitem{kulhlivre}
Franz-Viktor Kuhlmann.
\newblock {\em \textit{Valuation Theory}}.
\newblock \textup{http://math.usask.ca/~fvk/Fvkbook.htm}, 2011.

\bibitem{jcthese}
Jean-Christophe San~Saturnino.
\newblock {\em \textit{Th\'eor\`eme de Kaplansky effectif et uniformisation
  locale des sch\'emas quasi-excellents}}.
\newblock \textup{Th\`ese de Doctorat, Institut de Math\'ematiques de
  Toulouse}, 2013.

\bibitem{jccar0}
Jean-Christophe San~Saturnino.
\newblock \textit{Uniformisation locale des sch\'emas quasi-excellents de
  caract\'eristique nulle}.
\newblock {\em \textup{arXiv:1311.3525}}, 2013. Submitted.

\bibitem{jcdefaut}
Jean-Christophe San~Saturnino.
\newblock \textit{Defect of an extension, key polynomials and local
  uniformization}.
\newblock {\em arXiv:1412.7697}, 2014. Submitted.

\bibitem{vaquie3}
Michel Vaqui{\'e}.
\newblock \textit{Famille admissible de valuations et d\'efaut d'une
  extension}.
\newblock {\em \textup{J. Algebra}}, 311(2):859--876, 2007.

\bibitem{zarsam2}
Oscar Zariski and Pierre Samuel.
\newblock {\em Commutative Algebra II}.
\newblock Graduate Texts in Mathematics. Springer, 1976.

\end{thebibliography}
\end{document}